\documentclass[11pt]{article}
\usepackage[T1]{fontenc}
\usepackage{lmodern}
\usepackage[margin=1in]{geometry}
\usepackage{amsmath,amssymb,amsthm}
\usepackage{booktabs}
\usepackage{microtype}
\usepackage{graphicx}
\usepackage{float}
\usepackage{needspace}
\usepackage[hidelinks]{hyperref}

\makeatletter
\renewcommand\section{\@startsection{section}{1}{\z@}%
  {2.8ex plus -1ex minus -.2ex}{1.3ex plus .2ex}{\normalfont\Large\bfseries}}
\makeatother

\newtheorem{theorem}{Theorem}

\theoremstyle{definition}

\newcommand{\R}{\mathbb{R}}

\newcommand{\inner}[2]{\langle #1,#2\rangle}

\title{Improved lower bounds for the kissing number\\ in dimensions 49 to 61}
\author{Dmitry Kamenetsky}

\begin{document}
\maketitle

\begin{abstract}
We give explicit kissing configurations in $\R^d$ for every $d=49,\dots,61$. The lower bound
is improved over the best previously published value in all thirteen dimensions, by margins
from $232$ to $159{,}220$ points and by $506{,}461$ points in total. The configurations are
obtained from the antipodal block lift of Yang et al.\ \cite{Yang26} applied to the
$52{,}416{,}000$-point minimal shell of the extremal even unimodular lattice $P_{48p}$, using a class of $7{,}193$
minimal lines and, in each dimension, a kissing configuration of $\R^k$ ($k\le 15$) with the
largest available number of antipodal pairs. Two ingredients are new. The published class of
$7{,}077$ lines is grown to $7{,}193$ by local search in three stages, as a sequence of
exchanges between candidate lines and class lines; and the point count, which equals $52{,}416{,}000+2S+2t$ with $t$ the number of
antipodal pairs and $S$ the size of the \emph{union} of $t$ isometric images of the class, is
maximised as a maximum-coverage problem by greedy search with iterated local search over pools
of $4{,}000$ automorphic images, reaching $92$--$100\%$ of the trivial bound $7{,}193t$. All
results are verified in exact integer arithmetic, and the certificates are distributed with the
code. Code and data: \url{https://github.com/dimkadimon/KissingNumbers}.
\end{abstract}

\section{Introduction}

\emph{How many equal spheres can touch a given sphere, all of the same size and without
overlapping?} This question was put to Isaac Newton by David Gregory in 1694. Newton
held that the answer is $12$; Gregory suspected $13$. Newton was right, but a complete proof took more than three centuries \cite{Sch53}. The quantity is now called the \emph{kissing number} $K(d)$: the
largest number of non-overlapping unit spheres in $\R^d$ that can all touch one central unit
sphere, an arrangement of points on the unit sphere with pairwise distances at least one.

In the plane the answer is easy to see: six equal circles fit around a seventh, $K(2)=6$, and
the centres form a hexagonal pattern. In $3$D space $K(3)=12$, the centres sitting at the vertices of
a cuboctahedron, but here already the proof is delicate, and the visual intuition is misleading
--- the twelve neighbours leave visible gaps, which for a long time made it unclear whether a
thirteenth sphere could be squeezed in. Beyond three dimensions the problem becomes a
benchmark for the limits of both geometry and computation: despite much effort, $K(d)$ is known
exactly only in dimensions $1,2,3,4,8,24$ \cite{Mus08,Via17,CKMRV17}. In every other dimension
one has to settle for lower bounds, supplied by explicit constructions, and upper bounds. This
paper improves the lower bounds in dimensions $49$ to $61$.

Equivalently, one asks for the largest set of points on the unit sphere with pairwise distances
at least one, that is, a code with a fixed minimum separation. In the notation used below, $K(d)$
is the largest cardinality of a set $X\subset S^{d-1}$ with $\inner{x}{y}\le 1/2$ for distinct
$x,y$. The published tables of lower bounds traditionally stop at dimension $48$, where the extremal even unimodular lattice $P_{48p}$ of Nebe
\cite{Nebe98,Catalogue} has $52{,}416{,}000$ minimal vectors \cite{CS99}, and resume at
dimension $64$, where the code construction of Edel, Rains and Sloane \cite{ERS98} gives
$331{,}737{,}984$. The intermediate range $49$--$63$ was therefore, until recently, occupied
only by the trivial observation that a direct sum with a good low-dimensional configuration
gives $52{,}416{,}000+K(d-48)$.

That changed with two independent works. Yang et al.\ \cite{Yang26} introduced an
\emph{antipodal block lift} (their Theorem~2.5) which, applied to the minimal shell of
$P_{48p}$ together with a class of $7{,}077$ minimal lines, produces configurations in
dimensions $49$--$55$; and the repository \cite{KissNum} evaluated the same family in
dimensions $49$--$61$ with a class of $7{,}069$ lines and better tails. The two overlap in
$49$--$55$ and disagree in $56$--$61$, where \cite{KissNum} is the only published source.
Table~\ref{tab:main} lists the best previously published value in each dimension, with its
source.

This paper keeps the construction of \cite{Yang26} and changes one ingredient: the choice of
isometric images. Because each selected image contributes to the count through the size of the
\emph{union} of the images and not through its own size, the images must be chosen jointly. We
formalise this as a maximum-coverage problem, solve it with a greedy pass followed by exact
$1$-swap local search and iterated local search over pools of up to $4{,}000$ images, and
obtain an improvement in every dimension $49$--$61$.

A second lever points in the same direction. Every line of the class is worth
$2t$ points in dimension $48+k$, so enlarging the class by one line improves all thirteen
bounds at once. Section~\ref{sec:class} describes the searches that grow the published class of
$7{,}077$ lines to $7{,}193$ lines, in three stages; the images of
Section~\ref{sec:images} are then re-optimised for the enlarged class.

\section{The lift and the count}

We begin by recalling the construction and fixing notation. Let $G$ be the published integral
Gram matrix of $P_{48p}$, so that a lattice vector is an
integer row $c$ and its squared norm is $cGc^{T}$. The minimal shell has $52{,}416{,}000$
vectors of norm $6$ and is antipodal; a \emph{line} is a pair $\{c,-c\}$, represented by the
member whose first nonzero coordinate is positive.

A \emph{class} is a set $C_0$ of minimal lines with
\begin{equation}\label{eq:class}
|cGc'^{T}|\le 2 \quad\text{for distinct } c,c'\in C_0 ,
\end{equation}
which after normalisation is $\inner{c}{c'}\le 1/3$. The class used here has
$|C_0| = 7{,}193$ lines: it is obtained in Section~\ref{sec:class} from the $7{,}077$-line
class of \cite{Yang26} (one of whose elements is a $7{,}069$-line subclass used in
\cite{KissNum}) by adding $125$ lines and deleting $9$.

\begin{theorem}[Antipodal block lift; {\cite[Theorem 2.5]{Yang26}}]\label{thm:lift}
Let $C\subset S^{D-1}$ be a kissing configuration and $B_1,\dots,B_t\subset C$ pairwise
disjoint subsets with $\inner{x}{x'}\le 1/3$ for distinct $x,x'\in B_i$. Let $Y\subset S^{k-1}$
be a kissing configuration containing $t$ distinct antipodal pairs $\{\pm u_i\}$. Then
\[
X=\{(x,0):x\in C\setminus\textstyle\bigcup_i B_i\}
\;\cup\;
\{(\tfrac{\sqrt3}{2}x,\tfrac{\varepsilon}{2}u_i): x\in B_i,\ \varepsilon=\pm1\}
\;\cup\;
\{(0,y):y\in Y\}
\]
is a kissing configuration in $\R^{D+k}$ with $|X|=|C|+\sum_i|B_i|+|Y|$.
\end{theorem}

Apply this with $C$ the minimal shell of $P_{48p}$ and with blocks obtained from a class:
given isometries $A_1,\dots,A_t$ of the Gram matrix ($A_iGA_i^{T}=G$), assign each line of
$C_0A_i$ to the first image in which it occurs. The assigned parts are disjoint as line sets,
their union has
\begin{equation}\label{eq:S}
S=\Bigl|\bigcup_{i=1}^{t} C_0A_i\Bigr| \quad\text{lines}
\end{equation}
and every image is again a class by \eqref{eq:class}. Taking both signs of each assigned line
gives blocks $B_i$ of size $2|D_i|$ satisfying the hypothesis of Theorem~\ref{thm:lift}, and
\begin{equation}\label{eq:count}
K(48+k)\ \ge\ 52{,}416{,}000+2S+2t ,
\end{equation}
since a tail with $t$ antipodal pairs has $|Y|=2t$ points. This is equation (8.8) of
\cite{Yang26}.

Two consequences shape everything that follows. First, $t$ is determined only by the tail and
is independent of the lattice: any kissing configuration of $\R^k$ with $t$ antipodal pairs
will do, so the same class family serves every dimension. Second, the images enter
\eqref{eq:count} only through the single number $S$, which is at most $|C_0|t$ and is smaller
whenever two chosen images share a line. The class itself enters linearly: one more line is
worth two more points per image, that is $2t$ points in dimension $48+k$.

\section{Growing the class}\label{sec:class}

Fix the class problem: find as many minimal lines of $P_{48p}$ as possible with
\eqref{eq:class}. Two lines conflict when $|cGc'^T|>2$; a class is an independent set in the
conflict graph on the minimal lines, and the published class of $7{,}077$ lines is one such
set. Enlarging it is the only lever that scales every dimension at once, one added line being
worth $2t$ points in dimension $48+k$. We use three searches in sequence, each of which exchanges
candidate lines for class lines; the sizes grow from $7{,}077$ to $7{,}082$, then to $7{,}108$,
and finally to $7{,}193$.

\paragraph{Stage 1: conflict-filtered exchanges.} A line that conflicts with few class lines
can enter the class if those few conflicts leave. We draw a pool of $20$ random isometric images
of the class --- about $1.4\cdot10^{5}$ distinct lines --- and count, for every pool line, its
conflicts against the whole class. The distribution is sharply concentrated: over $141{,}140$
pool lines, $57$ lines conflict with exactly one class line, $105$ with two, $431$ with three,
$1{,}478$ with four, $3{,}619$ with five, and the median is ten. Only lines with at most five
conflicts are kept, which keeps the exchange problem small: if $R$ is the set of class lines hit
by the kept candidates and $A$ an independent set among them, then $(C_0\setminus R)\cup A$ is
again a class, and it is larger exactly when $|A|>|R|$. We solve this induced
maximum-independent-set instance by iterated local search --- delete $k\le6$ random lines of the
current solution, refill greedily --- and rebuild the pool around each improved class. Six
passes give $7{,}082$ lines.

\paragraph{Stage 2: streamed tightness search.} The second search works with the
\emph{tightness} of a candidate, its number of conflicts with the current class. Each batch
draws $20$ new automorphisms and forms all images of the current class lines; only candidates
with tightness at most $2$ are kept, which bounds memory while retaining the useful ones. On
that restricted candidate set the search maintains the tightness of every candidate
incrementally as lines enter and leave the solution.

Throughout, an \emph{$(a,b)$-swap} deletes $a$ class lines and inserts $b$ candidates in their
place, so the class grows by $b-a$ lines; an $(a,a)$-swap leaves the size unchanged. The search
repeatedly applies
\begin{itemize}
\item free adds, for candidates of tightness $0$;
\item $(1,2)$-swaps: two candidates of tightness $1$, conflicting with the same single class
      line, are compatible with each other, and replacing that line by the pair gains one;
\item random $(1,1)$-swaps on candidates of tightness $1$: these leave the size of the class
      fixed and let the search drift across a plateau of equally large classes.
\end{itemize}

\paragraph{Stage 3: $(r,r+1)$ exchanges and restarts.} The first two stages never consider a
candidate blocked by more than one class line. Stage 3 adds them: candidates that
share the \emph{same} set of $r\le3$ blocking class lines can be exchanged for that set whenever
$r+1$ of them are mutually compatible, which is a $(r,r+1)$-swap in the above sense; the blocker sets are
recomputed from scratch for each attempt, so the move does not depend on incremental counters.
The search is then restarted many times from the current best class with fresh automorphisms,
which turns it into a randomised hill climb with restarts.

\paragraph{Result.} Fifty batches of stage 2 give $7{,}108$ lines. The final class, after
stage 3, has $7{,}193$ lines, obtained from the published
$7{,}077$ by adding $125$ lines and deleting $9$. It passes the same exact checks --- norm
$cGc^{T}=6$, $|cGc'^{T}|\le 2$ for all $5\cdot10^{7}$ pairs, distinct lines --- and the extra
lines are worth up to $2t$ points each in dimension $48+k$: at $d=61$ alone the ceiling rises
from $60{,}589{,}782$ (with the published $7{,}077$ lines) to $60{,}717{,}876$.

None of the three searches is a certificate of maximality; all are hill climbs from a known
solution. The conflict statistics also explain why the naive alternative --- drop some lines at
random and refill --- fails: a random pool line conflicts with about ten class lines, so removing
a handful frees almost no candidate, while removing enough to help destroys the structured part
of the class faster than refilling restores it.

\section{Choosing the images}\label{sec:images}

For a fixed $t$, \eqref{eq:count} is maximised by maximising $S$. This is a maximum-coverage
instance: the ground set is the set of minimal lines of $P_{48p}$ that occur in at least one
image, the sets are the images $C_0A_i$, and exactly $t$ of them must be chosen.

\paragraph{Pool.} Distinct isometries are drawn as words of length $170$ in the three
published generators of $P_{48p}$, retaining those that satisfy $AGA^{T}=G$ exactly. Pools of
$400$--$4{,}000$ images are used; when re-optimising a previous family we seed the pool with
its images, which makes a pool of up to $4{,}577$. Each line of each image is replaced by a
sign-canonical $64$-bit key so that identical lines coincide; for $4{,}577$ images this is a
$32{,}000{,}000$-entry join, and the number of distinct lines in such a pool is about
$17.4$ million.

\paragraph{Greedy.} Images are added one at a time, each time the one covering the most
uncovered lines. This single pass produces the whole family: after $t$ additions we have a
valid configuration for dimension $48+k(t)$, so one run of the optimiser yields every
dimension at once.

\paragraph{Local search.} We use exact best-improvement $1$-swaps --- replace a chosen image by
an unchosen one and evaluate the true change in $|\bigcup_i C_0A_i|$ --- restricted to the $26$
most promising incoming and $8$ lowest-value outgoing images, followed by iterated local
search: remove a random sixth of the family, refill greedily, polish, and keep the result if
the union grew.

\paragraph{What binds.} The pool binds, not the local search. On the $7{,}193$-line class the
search reaches $S=3{,}820{,}380$ at $t=577$, i.e.\ $92.05\%$ of $|C_0|t$; the small-$t$ entries
sit at $99.95\%$ or above, where $1$-swaps alone close everything that is left. Starting the
search from an existing family and adding $4{,}000$ fresh images moves $S$ by thousands of
lines, whereas local search on a fixed pool moves it by tens.

\paragraph{Ablation.} Table~\ref{tab:ablation} separates the two levers: the families are held
fixed and only the class changes. Every stage of the growth from $7{,}077$ to $7{,}193$ lines
increases $S$ at every $t$, and at $t=577$ the last stage alone is worth $42{,}202$ lines of
union, out of the $85\cdot577=49{,}045$ that its $85$ new class lines could supply in the
absence of overlap; the shortfall is the overlap the larger class forces, not a failure of the
search.

\begin{table}[H]
\centering
\small
\begin{tabular}{@{}rrrrrrrr@{}}
\toprule
$d$ & $t$ & $|C_0|=7{,}077$ & $7{,}082$ & $7{,}108$ & $7{,}193$ &
      $S/(7{,}077t)$ & $S/(7{,}193t)$ \\
\midrule
49 & 1 & 7{,}077 & 7{,}082 & 7{,}108 & 7{,}193 & $100.000\%$ & $100.000\%$ \\
50 & 3 & 21{,}231 & 21{,}246 & 21{,}324 & 21{,}579 & $100.000\%$ & $100.000\%$ \\
51 & 6 & 42{,}460 & 42{,}490 & 42{,}646 & 43{,}156 & $99.995\%$ & $99.995\%$ \\
52 & 12 & 84{,}886 & 84{,}946 & 85{,}257 & 86{,}276 & $99.955\%$ & $99.954\%$ \\
53 & 20 & 141{,}360 & 141{,}459 & 141{,}977 & 143{,}672 & $99.873\%$ & $99.869\%$ \\
54 & 36 & 253{,}981 & 254{,}159 & 255{,}082 & 258{,}119 & $99.690\%$ & $99.680\%$ \\
55 & 63 & 442{,}905 & 443{,}214 & 444{,}819 & 450{,}091 & $99.339\%$ & $99.323\%$ \\
56 & 120 & 837{,}045 & 837{,}625 & 840{,}646 & 850{,}575 & $98.564\%$ & $98.542\%$ \\
57 & 153 & 1{,}062{,}166 & 1{,}062{,}901 & 1{,}066{,}697 & 1{,}079{,}243 & $98.096\%$ & $98.066\%$ \\
58 & 255 & 1{,}744{,}279 & 1{,}745{,}473 & 1{,}751{,}636 & 1{,}771{,}982 & $96.656\%$ & $96.607\%$ \\
59 & 302 & 2{,}051{,}070 & 2{,}052{,}463 & 2{,}059{,}649 & 2{,}083{,}435 & $95.968\%$ & $95.910\%$ \\
60 & 378 & 2{,}538{,}709 & 2{,}540{,}409 & 2{,}549{,}175 & 2{,}578{,}309 & $94.901\%$ & $94.827\%$ \\
61 & 577 & 3{,}763{,}184 & 3{,}765{,}682 & 3{,}778{,}178 & 3{,}820{,}380 & $92.157\%$ & $92.049\%$ \\
\bottomrule
\end{tabular}
\caption{The ablation of Section~\ref{sec:images}. The families are those of Table~\ref{tab:main}
and only the class changes --- the class is grown by $+5$, $+26$ and $+85$ lines in the three
stages, so each column shows the effect of the class alone on the union of a fixed set of
images. Every stage increases $S$ at every $t$; at $t=577$ the last stage adds $42{,}202$ lines
of union out of the $49{,}045$ that its $85$ new class lines could supply, the difference being
the overlap that a larger class forces. The last two columns give the fraction of $|C_0|t$
attained with the smallest and with the largest class: the small decrease as the class grows is
that same overlap, which is why the marginal line pays slightly less than $2t$ points.}
\label{tab:ablation}
\end{table}

\paragraph{Why the images cannot simply be taken disjoint.} One might hope for a family whose
images are pairwise disjoint, which would attain $S=|C_0|t$ exactly. Write
$o(g)=|C_0\cap C_0g|$ for an automorphism $g$; the loss $|C_0|t-S$ is
$\sum_{i<j}o(A_iA_j^{-1})$ at first order. For $120$ random automorphisms the mean of $o(g)$ is
$2.40$ (histogram: $0$ in $12$ cases, $1$ in $21$, $2$ in $37$, $3$ in $26$, $4$ in $12$,
$5$ in $7$, $6$ in $3$, $7$ in $2$), so a random family loses about
$2.4\binom{t}{2}$ lines, which is the same order as the losses we observe. Cyclic families
behave differently: if $g$ has odd order $m$, the union of the $m$ powers reaches
$99.6$--$99.9\%$ of $m\cdot 7{,}077$ (for $m=11$, $121$ lines lost out of $77{,}847$), but
every even-order element we sampled had $g^{m/2}=-\mathrm{Id}$, so its cyclic family contains
only $m/2$ distinct images and the union collapses to at most half the ceiling. Disjointness by
construction therefore needs a family with no short relations, not a single generator; this is
the same obstruction that the ceiling analysis of \cite{KissNum} describes, and we have not
overcome it.

\needspace{15\baselineskip}
\section{Tails}

By Theorem~\ref{thm:lift} the tail may be any kissing configuration of $\R^k$ containing $t$
antipodal pairs; it does not have to be a lattice shell. Table~\ref{tab:tails} lists the tails
used, with the number of points $|Y|=2t$ and of antipodal pairs $t$. All of them are antipodal,
which is what the lift requires: it consumes antipodal pairs, not points.

\begin{table}[H]
\centering
\small
\setlength{\tabcolsep}{4pt}
\begin{tabular}{@{}lccccccccccc@{}}
\toprule
$k=d-48$ & 3 & 4 & 5 & 6 & 7 & 8 & 9 & 10 & 11 & 12 & 13 \\
\midrule
tail & $A_3$ & $D_4$ & $D_5$ & $E_6$ & $E_7$ & $E_8$ & $P_{9a}$ & $P_{10b}$ & $P_{11c}$ & $K_{12}$ & $P_{13}$ \\
points $|Y|=2t$ & 12 & 24 & 40 & 72 & 126 & 240 & 306 & 510 & 604 & 756 & 1154 \\
pairs $t$ & 6 & 12 & 20 & 36 & 63 & 120 & 153 & 255 & 302 & 378 & 577 \\
\bottomrule
\end{tabular}
\caption{The tails used, by dimension $k=d-48$ of the added space. $A_1,A_2$ (not listed) give
$t=1,3$; the root systems $A_3$--$E_8$ give the classical values; $P_{9a},P_{10b},P_{11c},P_{13}$
are the record kissing configurations of Cohn \cite{Cohn}; $K_{12}$ is the $756$ minimal vectors
of the Coxeter--Todd lattice. Every one of them is antipodally symmetric, as the lift requires,
which fails for the $841$-point record configuration of $\R^{12}$, so the lattice value $378$ is
used there instead. Each direction system was re-verified here in exact arithmetic: constant
norm, antipodal, all pairwise inner products $\le1/2$, numerical rank $k$.}
\label{tab:tails}
\end{table}

\section{Results}

\begin{table}[H]
\centering
\small
\begin{tabular}{@{}rrrrrrr@{}}
\toprule
$d$ & $k$ & $t$ & previous best & source & \textbf{this work} & margin \\
\midrule
49 & 1 & 1 & 52{,}430{,}156 & \cite{Yang26} & \textbf{52{,}430{,}388} & $232$ \\
50 & 2 & 3 & 52{,}458{,}468 & \cite{Yang26} & \textbf{52{,}459{,}164} & $696$ \\
51 & 3 & 6 & 52{,}500{,}930 & \cite{Yang26} & \textbf{52{,}502{,}324} & $1{,}394$ \\
52 & 4 & 12 & 52{,}585{,}772 & \cite{Yang26} & \textbf{52{,}588{,}576} & $2{,}804$ \\
53 & 5 & 20 & 52{,}698{,}696 & \cite{Yang26} & \textbf{52{,}703{,}384} & $4{,}688$ \\
54 & 6 & 36 & 52{,}923{,}774 & \cite{Yang26} & \textbf{52{,}932{,}310} & $8{,}536$ \\
55 & 7 & 63 & 53{,}301{,}140 & \cite{Yang26} & \textbf{53{,}316{,}308} & $15{,}168$ \\
56 & 8 & 120 & 54{,}085{,}392 & \cite{KissNum} & \textbf{54{,}117{,}390} & $31{,}998$ \\
57 & 9 & 153 & 54{,}534{,}176 & \cite{KissNum} & \textbf{54{,}574{,}792} & $40{,}616$ \\
58 & 10 & 255 & 55{,}893{,}422 & \cite{KissNum} & \textbf{55{,}960{,}474} & $67{,}052$ \\
59 & 11 & 302 & 56{,}505{,}668 & \cite{KissNum} & \textbf{56{,}583{,}474} & $77{,}806$ \\
60 & 12 & 378 & 57{,}477{,}123 & \cite{KissNum} & \textbf{57{,}573{,}374} & $96{,}251$ \\
61 & 13 & 577 & 59{,}898{,}694 & \cite{KissNum} & \textbf{60{,}057{,}914} & $159{,}220$ \\
\bottomrule
\end{tabular}
\caption{Lower bounds for $K(d)$, $49\le d\le 61$. The margin is the difference from the best
previously published value; total margin $506{,}461$. In dimensions $49$ and $50$ the images happen to be
pairwise disjoint ($S=7{,}193t$), so the ceiling of the construction is attained exactly.}
\label{tab:main}
\end{table}

Table~\ref{tab:main} gives the results. In each dimension the value is
$52{,}416{,}000+2S+2t$ with the $S$ and $t$ of Table~\ref{tab:entries}.

\begin{table}[!ht]
\centering
\small
\begin{tabular}{@{}rrrrrr@{}}
\toprule
$d$ & $t$ & $S$ & max $7{,}193t$ & loss & $S/(7{,}193t)$ \\
\midrule
49 & 1 & 7{,}193 & 7{,}193 & 0 & $100.000\%$ \\
50 & 3 & 21{,}579 & 21{,}579 & 0 & $100.000\%$ \\
51 & 6 & 43{,}156 & 43{,}158 & 2 & $99.995\%$ \\
52 & 12 & 86{,}276 & 86{,}316 & 40 & $99.954\%$ \\
53 & 20 & 143{,}672 & 143{,}860 & 188 & $99.869\%$ \\
54 & 36 & 258{,}119 & 258{,}948 & 829 & $99.680\%$ \\
55 & 63 & 450{,}091 & 453{,}159 & 3{,}068 & $99.323\%$ \\
56 & 120 & 850{,}575 & 863{,}160 & 12{,}585 & $98.542\%$ \\
57 & 153 & 1{,}079{,}243 & 1{,}100{,}529 & 21{,}286 & $98.066\%$ \\
58 & 255 & 1{,}771{,}982 & 1{,}834{,}215 & 62{,}233 & $96.607\%$ \\
59 & 302 & 2{,}083{,}435 & 2{,}172{,}286 & 88{,}851 & $95.910\%$ \\
60 & 378 & 2{,}578{,}309 & 2{,}718{,}954 & 140{,}645 & $94.827\%$ \\
61 & 577 & 3{,}820{,}380 & 4{,}150{,}361 & 329{,}981 & $92.049\%$ \\
\bottomrule
\end{tabular}
\caption{The union sizes behind Table~\ref{tab:main}. The \emph{loss} is the number of lines
counted twice: the $t$ images hold $7{,}193t$ lines in total, but a line occurring in two or more
of them is one line of their union, so $S=7{,}193t-\text{loss}$, and each lost line costs two
points. The last column is the fraction of the trivial maximum $7{,}193t$ attained. The loss grows
quadratically in $t$, as expected from $\sum_{i<j}o(A_iA_j^{-1})$ with $o(g)\approx2.4$ for a
random automorphism $g$.}
\label{tab:entries}
\end{table}

The margins are largest where the published families were furthest from their own ceiling. For
$d=61$ the value in \cite{KissNum} corresponds to $S=3{,}740{,}770$, i.e.\ $91.6\%$ of
$7{,}077t$, while the family given here reaches $92.05\%$ of $7{,}193t$ and gains $79{,}610$
lines --- worth $159{,}220$ points, since every line is worth two points. At most
$125\cdot2\cdot577 = 144{,}250$ of those can be attributed to the new class lines, so at least
$14{,}970$ come from the image selection alone. The split between the two levers depends on $t$: at $t=1,3$ the family sits exactly at the
ceiling, so the whole margin is the class growth ($116\cdot2t$); for $t\ge6$ the measured margin
already exceeds $116\cdot2t$, the excess coming from the image selection. In the low dimensions the gains
are accordingly small in absolute terms: at $d=51$ the ceiling allows $t\cdot7{,}193=43{,}158$
lines and only two of them are lost.

\begin{figure}[t]
\centering
\includegraphics[width=\textwidth]{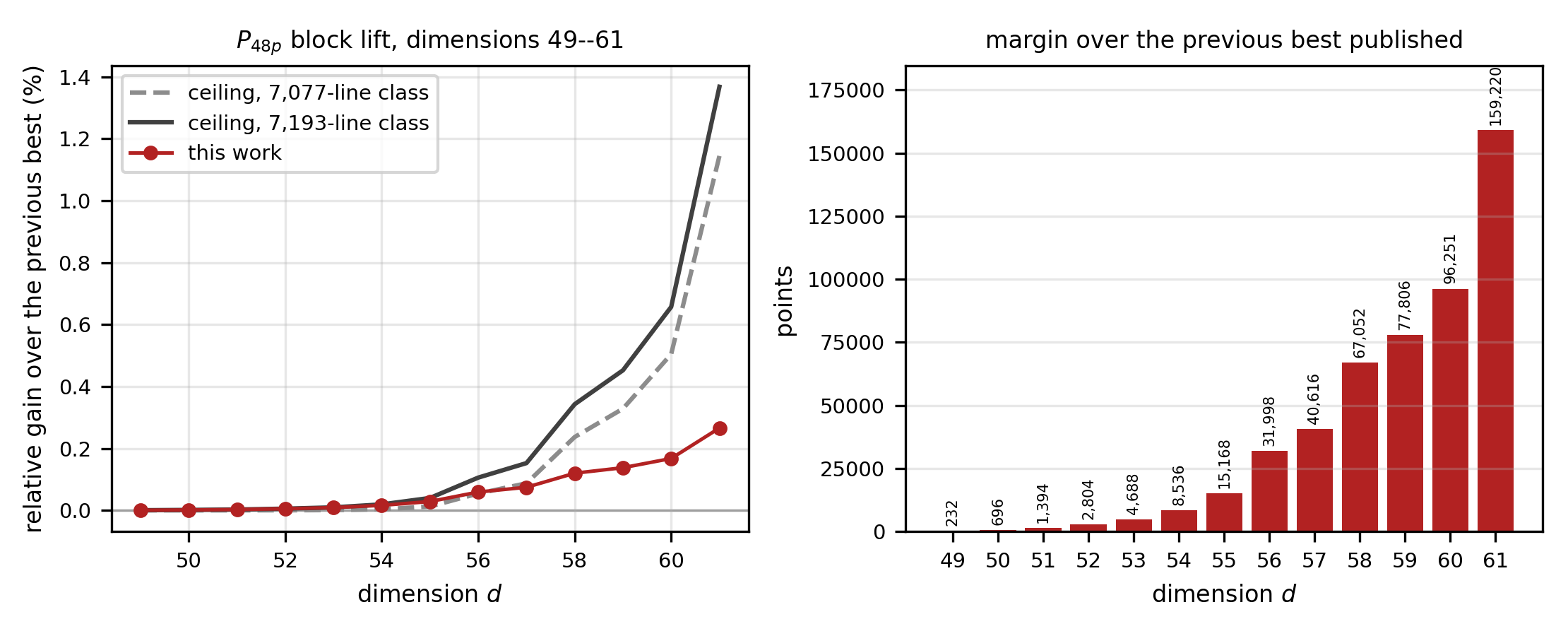}
\caption{Left: the margin of Table~\ref{tab:main} over the best previously published value, as a
percentage of that value (the line at $0\%$ is the previous best), together with the same
quantity for the ceilings $52{,}416{,}000+2\,|C_0|\,t+2t$ of \eqref{eq:ceiling} with
$|C_0|=7{,}077$ and $|C_0|=7{,}193$.
The curves are far below their ceilings, and the gap widens with $d$, which is where the remaining
room lies. Right: the margin in points; total $506{,}461$.}
\label{fig:main}
\end{figure}

Two independent checks confirm that these are real gains rather than bookkeeping. First, the
ceiling of the construction is exact and known: with every one of the $t$ classes at the maximum
$7{,}193$ lines the totals would be $52{,}416{,}000+2\cdot7{,}193t+2t$, and no family can
exceed them; every entry of Table~\ref{tab:main} is at or below its ceiling, and the two
entries at the ceiling are attained exactly. The ablation of Table~\ref{tab:ablation} separates
the two levers: its oldest-class column, with the published $7{,}077$ lines and the same families,
already reaches $3{,}763{,}184$ lines at $t=577$, against the $3{,}740{,}770$ of \cite{KissNum};
crediting all eight extra class lines to the class still leaves at least $17{,}798$ lines, that
is $35{,}596$ points, due to the choice of images. Second, a random control --- the same number of
images drawn at random rather than chosen by coverage --- loses roughly $2.4\binom{t}{2}$
lines, since the mean of $o(g)$ over random automorphisms is $2.40$; at $t=577$ that is of the
order of $4\cdot10^{5}$ lines, more than the $329{,}981$ lost by the optimised family. The
coverage optimisation, and the class growth behind it, are doing real work.
Figure~\ref{fig:main} shows the same margins in points, and relative to the previous best and to
the ceilings of \eqref{eq:ceiling}; Section~\ref{sec:stop} quantifies the room left above them.

\section{Where this method stops}\label{sec:stop}

The construction is capped by
\begin{equation}\label{eq:ceiling}
K(48+k)\ \le\ 52{,}416{,}000+2\,|C_0|\,t+2t ,
\end{equation}
which for $|C_0|=7{,}193$ and $d=61$ ($t=577$) is $60{,}717{,}876$. Our value is $98.9\%$ of
that wall, so at most $659{,}962$ further points are available in this dimension from any image
selection. Note that the wall moves with the class: the ceiling quoted in \cite{KissNum} for
$d=61$ is $60{,}574{,}780$, which is \eqref{eq:ceiling} with their $7{,}069$-line class. (With the
$7{,}077$-line class of \cite{Yang26} it would be $60{,}589{,}782$.) Every line found
by the searches of Section~\ref{sec:class} raises the ceiling by $2t$; the $116$ net lines we add are worth up to
$116\cdot2\cdot577 = 133{,}864$ points there, and a class of $8{,}250$ lines would put
$d=61$ above $61.9$ million. This is the ``class problem'', open in the same sense for the Leech lattice and
for $\Gamma_{72}$ \cite{KissNum}.

\paragraph{Future work.} Three directions seem most promising. First, the class: the searches of
Section~\ref{sec:class} were still gaining lines when they were stopped, and nothing in them
addresses the class problem itself, for which an upper bound would be as valuable as a larger
class; comparison with a spectral or Caro--Wei-type bound (which guarantees only $712$ lines for
$P_{48p}$ \cite{KissNum}) would show how far a local search can plausibly go. Second, the images:
the loss $|C_0|t-S$ grows quadratically in $t$ and eats $7.95\%$ of the ceiling at $d=61$; a
family with few short relations, which would make the images nearly disjoint, is the natural
target, and we have only ruled out the cyclic families that came to mind. Third, the framework:
the same lift applies to any lattice shell with a large class of lines at inner product
$\le1/3$, so the Leech lattice and $\Gamma_{72}$ \cite{KissNum} are the obvious next test beds,
and the tail side is limited only by the record antipodal kissing configurations of small
dimensions.

Above $d=61$ the method is not competitive and we make no claim there. Dimensions $62$ and $63$
are dominated by the Edel--Rains--Sloane code construction \cite{ERS98}, which gives
$71{,}310{,}732$ and $138{,}419{,}844$; \eqref{eq:ceiling} with the largest antipodal tails
there, $t=966$ and $t=1{,}282$ \cite{Cohn}, gives
only $66{,}314{,}808$ and $70{,}861{,}416$, below even the level-$0$ sign code of that chain. In
dimension $64$ and beyond \cite{ERS98} gives $331{,}737{,}984$ and more, which is five times
what this family can reach.

\section{Verification}

Each configuration is specified by the $t$ integer matrices $A_1,\dots,A_t$ and is checked from
scratch by \texttt{code/verify.py}, which
\begin{enumerate}
\item checks $A_iGA_i^{T}=G$ for every $i$, in exact integer arithmetic without floating point;
\item applies \texttt{data/class\_patch.json} to the published class, which adds $125$ lines and
      deletes $9$, and checks the resulting class: $cGc^{T}=6$ for all $7{,}193$ rows and
      $|cGc'^{T}|\le 2$ for every pair of distinct rows, and that the $7{,}193$ lines are
      distinct;
\item checks that no two selected images coincide as line sets;
\item recomputes $S=|\bigcup_i C_0A_i|$ \emph{without hashing}: the sign-canonical integer rows
      of all $t$ images are stacked and de-duplicated with \texttt{np.unique(axis=0)}, so the
      union size is exact by construction and independent of the search code;
\item reports $N=52{,}416{,}000+2S+2t$.
\end{enumerate}
The script \texttt{code/certificate.py} runs the same checks over the whole table in one pass and
writes \texttt{results/certificate.txt}, one line per dimension, ending with an assertion that
every entry agrees with \texttt{results/verified.json}, and records the SHA-256 checksums of
\texttt{data/class\_patch.json} and \texttt{results/verified.json}, so that a revision can be
cited unambiguously; pool size, word length and random seed behind each family are recorded in
\texttt{results/pool\_stats.json}.
The premises taken from the literature are only two: that the minimal shell of $P_{48p}$ has
$52{,}416{,}000$ vectors of norm $6$ \cite{CS99,Catalogue}, and Theorem~\ref{thm:lift}
\cite{Yang26}. The $7{,}077$-line class is the one supplied with \cite{Yang26}, and the added
lines are those of \texttt{data/class\_patch.json}; the property \eqref{eq:class} of the
resulting class is re-verified here rather than assumed.

\section*{Contribution statement}

The author proposed the problem, the overall approach and the specific experiments, and
reviewed all results. The research code, the search and the numerical experiments reported in
this paper were written and executed by an AI agent (Arena.ai Agent Mode), including the
implementation of the maximum-coverage search, of both class searches, of the exact verification
program, the comparison with the published bounds, and the preparation of this manuscript. The author checked the
verification code and the distributed certificates.

\end{document}